\documentclass[12pt]{article}
\usepackage{graphics}
\usepackage{graphicx}
\usepackage{amsmath,amsxtra,amssymb,latexsym, amscd}
\usepackage[left=3cm,right=2.5cm,top=2.5cm,bottom=2.5cm]{geometry}
\newtheorem{thm}{Theorem}
\newtheorem{cor}[thm]{Corollary}

\newtheorem{lem}[thm]{Lemma}
\newtheorem{prop}[thm]{Proposition}

\DeclareMathOperator{\ra}{rank}
\DeclareMathOperator{\diag}{diag}
\DeclareMathOperator{\card}{card}
\begin{document}
\baselineskip=17pt
\title{\bf Positive definite preserving linear transformations on symmetric matrix spaces\\
 }

\author{\bf Huynh Dinh Tuan-Tran Thi Nha Trang-Doan The Hieu\thanks{The authors are supported in part by a Nafosted Grant.}\\
\\
\sl  Hue Geometry Group\\
\sl  College of Education, Hue University\\
 34 Le Loi, Hue, Vietnam\\
\sl  dthehieu@yahoo.com
}
\maketitle

\begin{abstract}
Base on some simple facts of Hadamard product, characterizations of positive definite preserving linear transformations on real symmetric matrix spaces with an additional assumption ``$\ra T(E_{ii})=1, i=1,2,\ldots, n$'' or ``$T(A)>0\Rightarrow A> 0$'', were given.
\end{abstract}
\noindent {\bf AMS Subject Classification (2000):}  {Primary 15A86; Secondary 15A18, 15A04.}\\
{\bf Keywords:} {Linear preserver problems, Hadamard product, Symmetric matrix, Positive definite.}
\vskip 1cm

%==================================================
\section{Introduction}
In the recent years, one of active topics in the matrix theory is the linear preserver problems (LPPs). These problems  involve linear transformations on matrix space that have special properties, leaving some functions, subsets, relations \ldots invariant.
Many LPPs were treated while a lot of another ones are still open. For more details about LPPs: the history, the results and open problems we refer the reader to \cite{horn}, \cite{chi_tsi}, \cite{Li_Pi}, \cite{Loe}, and references therein.

On real symmetric or complex Hermitian matrices, one can consider the LPP of inertia. We say $A\in G(r,s,t)$ if $r, s, t$ are the numbers of positive, negative and zero eigenvalues of a matrix A, respectively.  The LPP of inertia asks to characterize a linear transformation $T$ preserving $G(r,s,t),$ i.e.
$$T(G(r,s,t))\subset G(r,s,t).$$

It is well-known (see \cite{Loe}) that a linear preserver of $G(r,s,t)$ on Hermitian matrices are of the form $A\mapsto WAW^*$ or $A\mapsto WA^tW^*$ for some invertible $W$ unless the cases of $rs=0$ or $r=s.$

The preservers of $G(n,0,0),$ the class of positive definite matrices, on complex Hermitian or real symmetric matrices are not known although one can easily find many such non-standard ones,  for example (on complex Hermitian matrices)
$$A\mapsto\sum_{i=1}^rW_iAW_i^*$$
where $W_i, i=1,2,\ldots, r$ are invertible matrices. In general, these can not be reduced to a single congruence  while there are some that can not be a sum of congruences (see \cite[Section 5]{Li_Pi}).

In this paper, we consider the LPP of $G(n,0,0)$ on real symmetric matrices. For simplicity, sometimes we write $A>0\ (A\ge 0)$ instead of ``$A$ is positive definite'' (``positive semi-definite'') and denote by $S_n(\Bbb R)$ the space of real symmetric matrices of size $n.$ Theorem \ref{thmmain} characterizes the linear transformations preserving both $G(n,0,0)$ and $\ra E_{ii},\ i=1,2,\ldots, n.$ They are of the form $A\mapsto W(H\circ A)W^t$ where $W$ is an invertible matrix, $H$ is a positive semi-definite matrix,\ $\diag H=\diag I_n,$  and the symbol ``$\circ$'' stands for the Hadamard product. As far as we known, this is  the first time, the Hadamard  product appears in a linear preserver problem.

If a preserver on real symmetric matrices $T$ is of the form $A\mapsto WAW^t,$ where $W$ is invertible, then it is obviously $T(A)>0$ implies $A> 0.$ The Theorem \ref{thmmain2} shows that, this condition is also a sufficient one.
Indeed, Theorem \ref{thmmain2} characterizes the ones preserving $G(n,0,0)$ with additional condition ``$T(A)>0 \Rightarrow A> 0$''. It is proved that this condition  implies $\ra T(E_{ii})=1,\ i=1,2,\ldots, n$ and hence $T$ has the form as in Theorem \ref{thmmain}. But in this case, the matrix $H$ is of rank 1 and therefore we can express $T$ in the standard form $A\mapsto WAW^t.$
The proof of the theorem also shows that the condition ``$T(A)>0 \Rightarrow A> 0$''
 is equivalent to ``$\det T(A)\ne 0 \Rightarrow \det A\ne 0$'',
 i.e. ``preserving singularity'';  ``$T({\cal A}) = {\cal A}$'' or ``preserving the set of all singular, positive semi-definite symmetric matrices''.

 Here ${\cal A}$ is the set of all positive definite symmetric matrices.

%=============================================
\section{Some useful lemmas}

It is well-known that a symmetric matrix $A$ is diagonalizable and associated with a quadratic form $q({\bf x})=\langle A{\bf x},{\bf x}\rangle,\ {\bf x}=(x_1, x_2, \ldots, x_n)\in \Bbb R^n.$ By a fundamental theorem of linear algebra, there exists an orthonormal basis $\{{\bf e_1}, {\bf e_2},\ldots, {\bf e_n}\}$ such that  $q$ can be brought to a diagonal form
$$q(x_1, x_2, \ldots, x_n)=\sum \lambda_i x_i^2.$$

In terms of matrices we can state

\begin{lem}\label{lem1}
For any $A\in S_{n}(\mathbb{R})$, there exists an orthogonal matrix $W$ such that $WAW^t$ is a diagonal matrix.
\end{lem}

If we set ${\bf x}_i=\sqrt{|\lambda_i|}{\bf e}_i,$ then we have two similar lemmas .

\begin{lem}\label{lem2}
For any $A\in S_{n}(\mathbb{R})$ of rank $r,$ there exists linear independent (pairwise orthogonal) vectors ${\bf x}_1, {\bf x}_2, \ldots, {\bf x}_r$ such that $A=\sum_{i=1}^{r}k_i{\bf x}_i{\bf x}_i^{t},\ k_i\in \{-1,1\}.$ Moreover if $A$ is positive semi-definite, then $A=\sum_{i=1}^{r}{\bf x}_i{\bf x}_i^t.$
\end{lem}

\begin{lem}\label{lem3}
For any positive semi-definite matrix $A\in S_{n}(\mathbb{R})$ of rank $r,$ there exists an invertible matrix $W$ such that $A=W(\sum_{i=1}^{r}E_{ii})W^t,$ where $E_{ii}={\bf e}_i{\bf e}_i^t.$
Moreover if $A$ is invertible then $A=WW^t.$
\end{lem}

Below are some more ones we need for the rest of the paper.

\begin{lem} \label{lemrank1}Let $A=\sum_{i=1}^{r}{\bf x}_i{\bf x}_i^t$ be a positive semi-definite matrix of rank $r$ in $S_n(\mathbb{R}).$ If every entry in its diagonal is non-zero, then there exists a vector \ ${\bf x}=(x_1, x_2, \ldots, x_n),\ x_i\ne 0,\  i=1,2,\ldots, n,$ such that $A={\bf x}{\bf x}^t+B,$ where $B$ is a positive semi-definite matrix of rank $r-1.$
\end{lem}
{\bf Proof.} Suppose that ${\bf x}_i=(x_{i1},x_{i2},\ldots, x_{in}).$ Let ${\bf u}_1=\alpha_1 {\bf x}_1+ \beta_1 {\bf x}_2=(u_{11}, u_{12}, \ldots, u_{1n})$ and  ${\bf v}_1=\beta_1 {\bf x}_1-\alpha_1 {\bf x}_2,$ where $\alpha_1$ and $\beta_1$ are chosen so that $\alpha_1, \beta_1>0;\ \alpha_1^2+\beta_1^2=1$ and whenever $x_{1k}\ne 0$ or $x_{2k}\ne 0, \ u_{1k}\ne 0.$ It is obviously that ${\bf x}_1{\bf x}_1^t+ {\bf x}_2{\bf x}_2^t= {\bf u}_1{\bf u}_1^t+{\bf v}_1{\bf v}_1^t.$  Let ${\bf u}_2=\alpha_2 {\bf u}_{1}+\beta_2 {\bf x}_{3}$ in such a way and so on, we get $A={\bf u}_{n-1}{\bf u}_{n-1}^t+\sum_{i=1}^{n-1}{\bf v}_i{\bf v}_i^t.$ Let ${\bf x}={\bf u}_{n-1}$ and $B=\sum_{i=1}^{n-1}{\bf v}_i{\bf v}_i^t,$ the lemma is proved.

\hfill$\Box$
\begin{lem} \label{lemrank2}Let $A=\sum_{i=1}^{r}k_i{\bf x}_i{\bf x}_i^t$ be a matrix of rank $r$ in $S_n(\mathbb{R}),\ k_i\in\{-1,1\},\ i=1,2,\ldots r$ and ${\bf x}_1,\cdots,{\bf x}_r$ are linear independent. If ${\bf x}\in\langle {\bf x}_1,\cdots,{\bf x}_r\rangle,$ where $\langle {\bf x}_1,\cdots,{\bf x}_r\rangle$ is the linear subspace generated by ${\bf x}_1,\cdots,{\bf x}_r,$ then $\ra (A+{\bf x}{\bf x}^t)\le \ra A.$
\end{lem}
{\bf Proof.}
It is clear that for any ${\bf u}\in\langle {\bf x}_1, {\bf x}_2, \ldots, {\bf x}_r\rangle^\perp,\ (A+{\bf x}{\bf x}^t)({\bf u})=0.$ This means that $\ra (A+{\bf x}{\bf x}^t)\le r= \ra A.$
\begin{lem} \label{lemA}
Let $A_1, A_2, \ldots, A_n$ be non-empty finite subsets of \ $\Bbb R^n,\ A_i\ne \{0\},\ i=1,2,\ldots, n.$ If for every $i\in\{1,2,\ldots, n\},$ $A_i\not\subset \sum_{j\ne i} A_j,$ then $\ra A_1= \ra A_2=\ldots=\ra A_n=1.$
\end{lem}

{\bf Proof.}
Suppose ${\bf v}_1, {\bf v}_2\in A_i$ are linear independent, then we can take a basic $B$ of $\sum A_j$ containing ${\bf v}_1, {\bf v}_2.$ Since $\card (B)\le n,$ there exists an index $k\ne i,$ such that $A_k$ is not belonging to $B,$ i.e. $A_k\subset\sum_{j\ne k}A_j$, a contradiction.
\hfill$\Box$

\begin{lem}\label{lemsemi}

A linear transformation (on symmetric matrices) preserving positive definiteness preserves  positive semi-definiteness.
\end{lem}

{\bf Proof.} Because $T$ is continuous and the topological closure of ${\cal A}$ is ${\cal B},$ where ${\cal B}$ is the set of all positive semi-definite symmetric matrices.
\hfill$\Box$

%==================================================

\section{The Hadamard product}
The Hadamard product is a simple matrix product, sometimes called the entrywise product. This product is much less widely understood although it has nice properties and some applications in statistics and physics (see \cite{horn}, \cite{john}).
For $(m\times n)$-matrices $A=(a_{ij})$ and $B=(b_{ij}),$ the Hadamard product of $A$ and $B$ is another $(m\times n)$-matrix, denoted by $A\circ B$ and defined by
$$  A\circ B=(a_{ij}b_{ij}).$$
It is easy to see that the Hadamard product is linear, commutative and has a nice relationship with diagonalizable matrices. Let $A=(a_{ij})$ be a diagonalizable matrix of size $n$ and $\lambda_1, \lambda_2, \ldots, \lambda_n$ be its eigenvalues. There exists an invertible matrix $W$ such that  $A=WDW^{-1},$ where $D=(d_{ij})$ is the diagonal matrix whose entries on the diagonal are $\lambda_i,$ i.e. $d_{ii}=\lambda_i,\ i=1,2,\ldots,n.$ We can verify the following
$$\begin{pmatrix}a_{11}\\ a_{22}\\ \vdots\\ a_{nn}\end{pmatrix}=(W\circ (W^{-1})^t\begin{pmatrix}\lambda_1\\ \lambda_2\\ \vdots\\ \lambda_n\end{pmatrix}.$$
\begin{thm}[The Schur product theorem]\ If $A, B$ are  positive semi-definite then so is $A\circ B.$
\end{thm}

Moreover, we have
\begin{prop}\label{prop0}
If $A\in S_n(\mathbb{R})$ is positive definite, $B\in S_n(\mathbb{R})$ is positive semi-definite and all entries on the diagonal of $B$ are non-zero then $A\circ B$ is positive definite.
\end{prop}
{\bf Proof.}
By Lemma \ref{lemrank1}, $B={\bf x}{\bf x}^t+C,$ where ${\bf x}=(x_1, x_2, \ldots, x_n),\ x_i\ne 0, i=1,2,\ldots, n$ and $C\ge 0.$ Suppose $A=\sum_{i=1}^{n}{\bf x}_i{\bf x}_i^t,$ where ${\bf x}_1,{\bf x}_2,\ldots, {\bf x}_n$ are linear independent. Because
$$\sum_{i=1}^{n}\alpha_i({\bf x}_i\circ {\bf x})=0\Leftrightarrow(\sum_{i=1}^{n}\alpha_i{\bf x}_i)\circ {\bf x}=0\Leftrightarrow\sum_{i=1}^{n}\alpha_i{\bf x}_i=0,$$
the vectors ${\bf x}_1\circ {\bf x},{\bf x}_2\circ {\bf x},\ldots,{\bf x}_n\circ {\bf x}$ are linear independent.

Thus,
$$\begin{aligned}
A\circ B&= A\circ({\bf x}{\bf x}^t)+A\circ C\\
&\ge A\circ({\bf x}{\bf x}^t)=\sum_{i=1}^{n}{\bf x}_i{\bf x}_i^t\circ({\bf x}{\bf x}^t)\\
&=\sum_{i=1}^{n}({\bf x}_i\circ{\bf x})({\bf x}_i\circ{\bf x})^t>0.\end{aligned}$$
\hfill$\Box$

The following proposition is useful for the proofs of the main theorems.

\begin{prop} \label{pro1}
Let ${\bf x}=(x_1, x_2, \ldots, x_n), {\bf y}=(y_1, y_2, \ldots, y_n)\in \Bbb R^n,$ where $x_i\ne 0, i=1,2,\ldots, n.$ If \ ${\bf x}, {\bf y}$ are linear independent, then there exists a singular positive semi-definite matrix $A$ of rank $n-1$ such that $A\circ ({\bf x}{\bf x}^t+{\bf y}{\bf y}^t)$ is positive definite.
\end{prop}
{\bf Proof.}
Without loss of generality, we can suppose that $\frac {y_1}{x_1}\ne \frac {y_n}{x_n}.$ Let ${\bf e}_1, {\bf e}_2, \ldots, {\bf e}_n$ be the standard basis of $\Bbb R^n$ and $A=\sum_{i=1}^{n-1}{\bf x}_i{\bf x}_i^t,$ where ${\bf x}_1={\bf e}_1+{\bf e}_n, {\bf x}_i={\bf e}_i, i=2,\ldots n-1.$ We can see that $A$ is a  positive semi-definte matrix of rank $n-1$ and $\{{\bf x}_1\circ {\bf y},{\bf x}_1\circ {\bf x},{\bf x}_2\circ {\bf x},\ldots, {\bf x}_{n-1}\circ {\bf x}\}$ is a basis of $\Bbb R^n.$ By Proposition \ref{prop0}, the proof is now proved because
$$\begin{aligned}
A\circ ({\bf x}{\bf x}^t+{\bf y}{\bf y}^t)&=(\sum_{i=1}^{n-1}{\bf x}_i{\bf x}_i^t)\circ({\bf x}{\bf x}^t+{\bf y}{\bf y}^t)\notag\\
&\geq(\sum_{i=1}^{n-1}{\bf x}_i{\bf x}_i^t)\circ({\bf x}{\bf x}^t)+({\bf x}_{1}{\bf x}_{1}^t)\circ({\bf y}{\bf y}^t)\notag\\
&=\sum_{i=1}^{n-1}({\bf x}_i\circ {\bf x})({\bf x}_i\circ {\bf x})^t+({\bf x}_{1}\circ {\bf y})({\bf x}_{1}\circ {\bf y})^t>0.\notag
\end{aligned} $$
\hfill$\Box$
%===================================
\section{The main theorems}

%====================================

\begin{thm}\label{thmmain}
Let $T: S_n(\mathbb{R})\longrightarrow S_n(\mathbb{R})$ is  a linear transformation preserving positive definiteness and  $\ra T(E_{ii})=1,\  i=1,2,\ldots, n.$ Then there exists an invertible matrix $W$ and a positive semi-definite matrix $H,\ \diag H=\diag I_n,$ such that for every $A\in S_n(\mathbb{R})$
\begin{equation}\label {main}
T(A)=W(H\circ A)W^t.\end{equation}
\end{thm}
{\bf Proof.}\
Since $T$ preserves positive definiteness, $T(I_n)$ is a positive definite matrix. Then there exists an invertible matrix $W_1$ such that $T(I_n)=W_1W_1^t.$
We can verify that the linear operator
$T_1(A)= W_1^{-1}T(A)(W_1^t)^{-1}$ also preserves positive definiteness, $\ra(T_1(E_{ii})=1,\ i=1,2,\ldots,n$ and moreover $T_1(I_n)=I_n$.\\
Suppose $T_1(E_{ii})={\bf u}_i{\bf u}_i^t,i=1,\ldots,n$. Since $I_n=T_1(I_n)=\sum T_1(E_{ii})=\sum {\bf u}_i{\bf u}_i^t>0;\ \{{\bf u}_1,{\bf u}_2, \ldots, {\bf u}_n\}$ is an orthonormal basis of $\Bbb R^n.$ Let $W_2=[{\bf u}_1,\cdots,{\bf u}_n],$ the orthogonal matrix whose $i$-th column is ${\bf u}_i$ and consider $T_2: S_n(\mathbb{R})\longrightarrow S_n(\mathbb{R})$, $T_2(A)= W_2T_1(A)W_2^t,\ \forall  A\in S_n(\mathbb{R})$. We can verify that $T_2$ has the same properties as  $T_1$'s  and moreover $T_2(E_{ii})=E_{ii},i=1,\ldots,n$.\\

 Let  $T_2(E_{ij}+E_{ji})=A.$ For every ${\bf y}=(y_1, y_2,\ldots, y_n),$ where $y_i=0$ (i.e. $\langle T_2(E_{ii}){\bf y},{\bf y}\rangle=\langle E_{ii}{\bf y},{\bf y}\rangle=0),$ we claim that
\begin{equation}\label{2}\langle A{\bf y},{\bf y}\rangle=0.\end{equation}

Indeed, suppose that there exists a such ${\bf y}\in\mathbb{R}^n$ such that $\langle A{\bf y},{\bf y}\rangle\ne 0$. We assume that $\langle A{\bf y},{\bf y}\rangle< 0.$ (The case of $\langle A{\bf y},{\bf y}\rangle> 0$ has a similar proof). Then we can choose  a small enough positive number $\epsilon$ and a big enough positive number $\beta$  (says $\beta\epsilon>1$) such that
$$ \langle A{\bf y},{\bf y}\rangle+\epsilon\langle(I-E_{ii}){\bf y},  {\bf y}\rangle <0,$$
and
$$X=(\beta-\epsilon)E_{ii}+\epsilon I +(E_{ij}+E_{ji})\ \ {\text {is positive definite}}.$$

But then $\langle T_2(X){\bf y},{\bf y}\rangle=  \langle A{\bf y}, {\bf y}\rangle+\epsilon \langle(I-E_{ii}){\bf y},{\bf y}\rangle<0,$ a contracdition.

The equality (\ref{2}) means that, all entries of $A$ are zeros unless ones lying on the $i$-th column or the $i$-th row. The equality (\ref{2}) also holds for ${\bf y}=(y_1, y_2,\ldots, y_n),$ where $y_j=0$ and therefore all entries of $A$ are zeros unless ones lying on the $j$-th column or the $j$-th row. Thus, for every $i,j,\ i\ne j$
  $$T_2(E_{ij}+E_{ji})=h_{ij}(E_{ij}+E_{ji}).$$

Set $H=(h_{ij})_{n\times n}, h_{11}=h_{22}\ldots=h_{nn}=1.$ Since $H=T({\bf 1}),$ where ${\bf 1}$ is the square matrix whose all entries are 1, $H$ positive semi-definite by Lemma \ref{lemsemi}. We can verify that
$$T_2(A)=A\circ H$$
for all $A\in S_n(\mathbb{R}).$
Let $W=W_2W_1,$ we have (\ref{main}).
It is easy to see that, if $T$ has the form (\ref{main}), then $T$ preserves positive definiteness, by Proposition \ref{prop0}, and $\ra T(E_{ii}) =1$ for all $i=1,2,\ldots, n.$
\hfill$\Box$
\vskip.5cm
%========================================================
%It is easy to see that, if $T(A)=WAW^t$ then
\begin{thm}\label{thmmain2}
 Let $T: S_n(\mathbb{R})\longrightarrow S_n(\mathbb{R})$ is a linear operator preserves positive definiteness. Then, $T$ satisfies the condition $T^{-1}({\cal A})\subset {\cal A}$ if and only if
 there exists an invertible  $W$ such that for every $A\in S_n(\mathbb{R})$

\begin{equation}\label{3} T(A)=WAW^t\end{equation}
\end{thm}

{\bf Proof.}
It is easy to see that, if $T(A)=WAW^t$ then $T^{-1}({\cal A})\subset {\cal A}.$

Now suppose that $T^{-1}({\cal A})\subset {\cal A}.$

First, we prove that $\ra T(E_{ii})=1,i=1,\cdots,n.$

If there exists an index $i$ such that $T(E_{ii})=\{0\},$ then $\sum_{j\not=i}T(E_{jj})=T(\sum_{j\not=i}E_{jj})=T(I_n)>0,$ a contradiction because $\det\sum_{j\not=i}E_{jj}=0.$ Thus, $T(E_{ii})\ne\{0\}, \ i=1,2,\ldots,n.$

Suppose $\ra T(E_{ii})=r_i\ge 1,i=1,\cdots,n.$ Since $T(E_{ii})$ is positive semi-definite (Lemma \ref{lemsemi}), $T(E_{ii})=\sum_{j=1}^{r_i}{\bf u}_{ij}{\bf u}_{ij}^t$. Denote by $A_i=\{{\bf u}_{ij}: j=1,2,\ldots, n\}.$

If there exists an index $i$ such that $A_i\subset \sum_{j\ne i}A_j,$ then by Lemma \ref{lemrank2}, $n=\ra T(I_n)\leq \ra\sum_{j\not=i}T(E_{jj}) ,$ a contradiction. Thus,  $A_i\not\subset \sum_{j\ne i}A_j$ and hence $\ra T(E_{ii})=1,\ i=1,\cdots,n$ by virtue of Lemma \ref{lemA}.

By Theorem \ref{thmmain}, there exists an invertible matrix $W_1$ such that
$$T(A)=W_1(H\circ A)W_1^t,\ \ \ \ \forall A\in S_n(\mathbb{R}).$$

Following the proof of Theorem \ref{thmmain}, $\diag H=\diag I_n$. Suppose that $\ra(H)=r,\ r>1$. By Lemma \ref{lemrank1} and because $\ra H>1,\ H={\bf x}{\bf x}^t+{\bf y}{\bf y}^t+B,$ where ${\bf x}, {\bf y}$ are linear independent, ${\bf x}=(x_1, x_2,\ldots,x_n),\ x_i\ne 0, i=1,2,\ldots n$ and $B$ is positive semi-definite.
By Proposition \ref{pro1}, there exists a singular (and positive semi-definite) matrix $A$ such that $A\circ ({\bf x}{\bf x}^t+{\bf y}{\bf y}^t)$ is positive definite. But $A\circ H\ge A\circ  ({\bf x}{\bf x}^t+{\bf y}{\bf y}^t)>0$ and hence $T(A)>0.$  This contradiction means that $\ra H=1$ and we have $H={\bf u}{\bf u}^t,$ where ${\bf u}=(u_1, u_2, \ldots, u_n),\ \ u_i\ne 0,\  i=1,2,\ldots, n.$ Since $\diag H=\diag I_n,\ u_i\in\{-1,1\}, i=1,2,\ldots, n.$ It is not hard to see that $A\circ H=(\sum u_iE_{ii})A (\sum u_iE_{ii})$ and set $W=W_1(\sum u_iE_{ii}).$ The theorem is proved.
\hfill$\Box$
\vskip.5cm

{\bf Remark.}
 Positive definite preservers of standard form (expressed by a single congruence) are only ones satisfying $T^{-1}({\cal A})\subset {\cal A}.$

Following the proof of Theorem \ref{thmmain2}, it is not hard to prove the following

\begin{cor}  The condition $T^{-1}({\cal A})\subset {\cal A}$ in the Theorem \ref{thmmain2} is equivalent to each one of the followings:
\begin{enumerate}
\item  $T^{-1}({\cal C})\subset {\cal C},$ where
    ${\cal C}$ is the set of all invertible symmetric matrices, i.e. $T$ preserves singularity.
\item $ T$ preserves the set of all singular, positive semi-definite symmetric matrices''.
\item $T({\cal A})={\cal A}.$
\end{enumerate}

\end{cor}

%%%%%%%%%%%%%%%%%%%%%%%%%%%
%{\bf Acknkowledgements.}


\begin{thebibliography}{99}

\bibitem {cao} Cao, C. and Tang, X., \textit{Determinant Preserving Transformations on Symmetric matrix Spaces}, Electronic Journal of Linear Algebra, Vol.11, 205-211 (2004).
\bibitem {horn} Horn, R. A. and Johnson, C. R., \textit{Topics in Matrix Analysis}, Cambridge University Press, Cambridge (1991).
\bibitem{john} Johnson, C., \textit{Matrix Theory and Applications}, American Mathematical Society, 1990.
\bibitem{chi_tsi} Li, C. K. and Tsing, N. K., \textit{Linear preserver problems: a brief introduction and some special techniques}, Directions in matrix theory (Auburn, AL, 1990). Linear Algebra Appl. 162/164  217-235 (1992).
\bibitem{Li_Pi} Li, C.-K. and Pierce, S., \textit{Linear preserver problems}, Amer. Math. Monthly, 108: 591-605 (2001).
\bibitem{Loe1} Loewy, R.,  \textit{Linear maps which preserve a balanced nonsingular inertia class}, Linear Algebra Appl. 134:165-179 (1990).
\bibitem{Loe2}Loewy, R.,  \textit{Linear maps which preserve an inertia class}, SIAM J. Matrix Anal. Appl. 11:107-112 (1990).
\bibitem{Loe} Loewy, R., \textit{A survey of linear preserver problems-chapter 3: Inertia preservers}, Linear Multilinear Algebra 33: 22–30 (1992).

%\bibitem{Pie} Pierce, S. et al., \textit{A survey of linear preserver problems}, Linear Multilinear Algebra 33: 1-129 (1992).
\bibitem{ste_Lei} Pierce, S. and Rodman, L., \textit{Linear Preservers of the class of Hermitian matrices with balanced inertia}, Siam J. Matrix Anal. Appl. Vol. 9, No. 4, 461-472 (1988).
\end{thebibliography}
\end{document}